\documentclass[12pt]{amsart}
\usepackage{verbatim} 
\usepackage{amsfonts,amssymb,amsmath}
\begin{document}
\pagestyle{plain}
\large
\newtheorem{thm}{Theorem}
\newtheorem{note}{Note}
\newtheorem{notation}{Notation}
\newtheorem{Def-Not}{Definition}
\newtheorem{prop}{Proposition}
\newtheorem{lemma}{Lemma}
\newtheorem{example}{Example}

\newtheorem{Definition}{Definition}[section]
\newtheorem{definition}{Definition}
\newtheorem{Proposition}{Proposition}[section]
\newtheorem{proposition}{Proposition}
\newtheorem{Remark}{Remark}[section]
\newtheorem{remark}{Remark}
\newtheorem{Restriction}{Restriction}[section]
\newtheorem{Corollary}{Corollary}[section]
\newtheorem{corollary}{Corollary}
\newtheorem{Theorem}{Theorem}[section]
\newtheorem{theorem}{Theorem}
\newtheorem{Lemma}{Lemma}[section]
\newcommand{\const}{\mbox{const.}}
\newtheorem{lem}{Lemma}[section]
\renewcommand{\thepage}{ \arabic{page}}

\def\CC{\mathbb{C}}
\def\NN{\mathbb{N}}
\def\CP1{\mathbb{C}{\mathbf P}^1}
\def\ZZ{\mathbb{Z}}
\def\RR{\mathbb{R}}
\def\OOO{\mathbf{O}}
\def\cL{{\mathcal  L}}
\def\gl1{\operatorname{gl}(1, \CC)}
\def\GL1{\operatorname{GL}(1, \CC)}
\def\SL2C{\operatorname{SL}(2,\mathbb{C})}
\def\GL2C{\operatorname{GL}(2,\mathbb{C})}
\def\sl2C{\operatorname{sl}(2,\mathbb{C})}
\def\glmC{\operatorname{gl}(m,\mathbb{C})}
\def\gl2C{\operatorname{gl}(2,\mathbb{C})}
\def\GLNC{\operatorname{GL}(N,\mathbb{C})}

\def\GLm{\operatorname{GL}(m)}
\def\SLmC{\operatorname{SL}(m,\mathbb{C})}
\def\slmC{\operatorname{sl}(m,\mathbb{C})}
\def\glN{\operatorname{gl}(N)}
\def\GLN{\operatorname{GL}(N)} 

\def\im{\operatorname{im}}
\def\dim{\operatorname{dim}}
\def\Res{\operatorname{Res}}
\def\Ad{\operatorname{ad}}
\def\tr{\operatorname{tr }}
\def\sn{\operatorname{sn }}
\def\cn{\operatorname{cn }}
\def\dn{\operatorname{dn }}
\def\const{ \mathfrak{const}}
\def\rank{ \operatorname{rank}}

\def\aA{\mathcal A}
\def\tildeV{\stackrel{\sim}{V}}
\def\tildeA{\stackrel{\sim}{\mathcal A}}

\title{
Rational version of Archimedes symplectomorphysm
and \\ birational Darboux coordinates on
coadjoint orbit of $\operatorname{GL}(N,{\mathbb C})$
}
\address{MPIM Bonn,  PDMI St.Petersburg}%
\email{mbabich@pdmi.ras.ru, misha.babich@gmail.com}%
\author{Mikhail V. Babich}

\begin{abstract}

A set of all linear transformations with a
fixed Jordan structure $J$ is
a symplectic manifold isomorphic to the coadjoint orbit $\mathcal O (J)$ of
$\operatorname{GL}(N,{\mathbb C})$.

Any linear transformation may be projected along its eigenspace
to (at least one) coordinate subspace of the complement dimension.
The Jordan structure $\tilde J$
of the image is defined by the Jordan structure $J$ of the pre-image,
consequently the projection
$\mathcal O (J)\to \mathcal O (\tilde J)$ is the mapping of the
symplectic manifolds.

It is proved that the fiber $\mathcal E$ of the projection
is a linear symplectic space  and the map $\mathcal O(J)
\stackrel{\sim}{\to} \mathcal E \times \mathcal O (\tilde J)$ is a
birational symplectomorphysm.

The iteration of the procedure gives the isomorphism between $\mathcal O (J)$
and the linear symplectic space, which is the direct product of all the fibers
of the projections.
The Darboux coordinates on $\mathcal O(J)$ are pull-backs of the canonical
coordinates on the linear spaces in question.

\end{abstract}

\maketitle
\tableofcontents

\setcounter{section}{-1}
\section{Introduction}

It was Archimedes who found that the ``proper'' coordinates for the element of the area of the sphere $\mathbf O({R})$ have a specific dual nature.

One coordinate is {\em the length}. It gives the position of the orthogonal projection of the parameterizing point on the diameter $\mathcal P$
connecting poles of the sphere.

The conjugated coordinate is {\em the angle}. This angle pa\-ra\-me\-tri\-zes the elements of 1-parametric subgroup $\mathcal Q\subset \operatorname{O}(3)$,
preserving the fibration $\mathbf O({R}) \to \mathcal P$ of the sphere on the circles $\mathcal C_{p}$ by the planes orthogonal to the diameter $\mathcal P\ni p$.

It is the cylindrical coordinates and the famous Archimedes area-preserving correspondence between the sphere $\mathbf O(R)$ and 
its circumscribing cylinder $\mathcal P\times \mathcal C_R$.
These sphere and cylinder were placed on the tomb of Archimedes at his request~\cite{Cicero}.

\vskip 3 ex
Let us demonstrate how the Archimedes method introduces (the standard) birational Darboux coordinates 
on the coadjoint orbit $\mathcal O$ of
$\operatorname{GL}(2, \CC)$.

\vskip 2 ex
We identify Lie algebra with its dual using a form $<A,B>= \tr AB$, and treat elements $A\in \mathcal O$ as matrices with the given Jordan structure. 
Let $\lambda_1, \lambda_2\in \CC$ be eigenvalues of $A\in \mathcal O$.
We put $\lambda_1=0$ and consider matrices with one zero eigenvalue, it is always possible to add a matrix proportional
to the unit matrix to the answer.

Let us consider any orbit $\mathcal O({J_R})\subset \operatorname{gl}(2, \CC)$, fixed by non-zero $J_R$ with the eigenvalues $0$ and $R$:

$$J_R=\left(\begin{array}{cc}
0&0\\
0&R
\end{array}
\right), \text{ if $R\ne 0$},\ \  J_0=\left(\begin{array}{cc}
0&1\\
0&0
\end{array}
\right).
$$
The orbit $\mathcal O({J_R})$ coincides with the (non-singular part of the) affine quadric
$$A\in {\mathcal O}({J_R}) \Leftrightarrow \det A=0, \tr A=R, A\ne 0.
$$
It is visual to consider the one-sheet hyperboloid $pY=X(R-X)$ in the 3-space of $X,Y,p$. The hyperboloid is fibrated on the parabolas $\mathcal C_p$:
$$\mathcal C_p: \ \ pY=X(R-X), \ \ p=\const, \ \ \ 
A= 
\left(\begin{array}{cc}
X&p\\
Y&R-X
\end{array}
\right)\in \mathcal O({J_R})$$

Consider the one-parametric subgroup $\mathcal Q\in\operatorname{GL}(N, \CC)$ preserving sections $p=\const$. It is a subgroup of unitriangular matrices
$\left(\begin{array}{cc}
1&0\\
q&1
\end{array}
\right)\in \mathcal Q$. 
Its elements shift the natural parameter on each parabola like 
the rotations of the Archimedes sphere shift the natural parameter (angle) of its sections (circles).

Variables $p,q$ parameterize a Zariski-open subset of $\mathcal O({J_R})$:
\begin{equation}\label{eq1-intro}
A(p,q)= 
\left(\begin{array}{cc}
1&0\\
q&1
\end{array}
\right)
\left(\begin{array}{cc}
0&p\\
0&R
\end{array}
\right)
\left(\begin{array}{cc}
1&0\\
-q&1
\end{array}
\right)=
\left(\begin{array}{cc}
-pq&p\\
-q(pq+R)&pq+R
\end{array}
\right)
\end{equation}

The calculation of the standard symplectic Lie-Poisson-Kirillov-Kostant form shows that $p,q$ are Darboux coordinates.

\vskip 2ex

In the present article we prove that the same trick in the combination with the simple iteration procedure gives the
birational Darboux coordinates on the (co)adjoined orbit $\mathcal O(J)$ corresponding to the matrix $J$ 
of any size and any Jordan structure.

\vskip 3ex
Parameterizations of coadjoint orbits have attracted attention of many authors. 
I want to mention the papers \cite{VN,Harnad,GHL,DubMaz,Krich,Fadd} that brought primary influence to bear on
the author. All these works were initiated by needs of the theory of integrable systems,
where coadjoint orbits are used for the construction of the phase spaces 
within the framework of the Hamiltonian formalism.

The papers \cite{GHL,Fadd} use the Gelfand-Zeitlin method. The explicit formulae are presented in \cite{GHL},
where  the authors introduce the nice parameterization of the generic orbit, where one family of the coordinates is formed by the eigenvalues
of the diagonal blocks of the matrix which is parameterized. The canonically conjugated coordinates 
can be easily calculated using vector-columns of the matrix and the corresponding eigenvectors of the blocks.
These coordinates are not rational but algebraic because it is necessary to find the eigenvalues of the matrices of all
sizes smaller than the initial one.

Another method was introduced in \cite{VN}. 
Really it is not method of the parameterization of $\operatorname{GL}(N,{\mathbb C})$ orbit, but some general 
scheme that can be applied to different problems.
The corresponding coordinates are called {\em the Spectral Darboux Coordinates}, 
see  \cite{DubMaz} where they are considered in detail. 
The method was applied to the parameterization of the $\operatorname{SL}(N,{\mathbb C})$-orbit in \cite{Harnad}.
It is not rational but algebraic too, it is necessary to solve algebraic equations of high order.

Consider the so called {\em Isomonodromic Coordinates} introduced in \cite{DubMaz}.
They are parameterize not an orbit but some manifold related to orbits again.
These coordinates are constructed for the orbits of the generic type and they are not rational but algebraic again.
Nevertheless the method of \cite{DubMaz} uses the cyclic process of the reduction of the
matrix equation of the first order to the scalar equation of the high order. 
This process is somewhat similar to the method developed in the present paper, 
but the connection between these two parameterizations is not understood yet,
it should be a subject for a future investigation.
\vskip 2ex

It is the rationality that is the fundamental property 
of the presented coordinates which is very important, at least for the investigation of isomonodromic deformations of systems 
of linear differential equations.

In this context the orbit is a Zariski-open domain of the phase-space of the corresponding algebraic Hamiltonian system.
The birationality of the transition functions is deeply connected with
the famous Painlev\'e-property of the isomonodromic deformations.
 
We consider not a generic but {\em the general} case of the structure of the orbit here. 
There are the orbits swept by degenerated matrices that are the cases of high importance 
because of the small dimensions of the corresponding orbits.
Such orbits can be treated as the phase spaces of ``more classical'', low-dimensional systems
immersed into ``more roomy'' 
high-dimensional spaces of matrices of higher sizes, that makes possible to find new approaches to old problems, see \cite{Boalch}.

\vskip 2ex
Let us turn to the subject of the paper. 
The crucial idea is the factorization of the matrix from the orbit on {\em the proper} 
triangular factors of a different nature. The idea belongs to  S.~E.~Derkachov and A.~N.~Manashov, 
they use it for the needs of quantum field theory \cite{DM}. 
Recently the method was applied for the parameterization of the orbits swept by the diagonalizable matrices  \cite{BD}. 
In the present article we are giving the evaluation of the method of \cite{DM,BD} to the general Jordan case.

\vskip 3ex

The main idea of the method can be demonstrated on the toy-example of $2\times 2$ case. 
Formula (\ref{eq1}) shows that, if we transform the first vector of the basis of 
the space to the eigenvector of the matrix $A$ by the transformation $\left(\begin{array}{cc}
1&0\\
q&1
\end{array}
\right)$, the Jordan form of the
resulting matrix $\left(\begin{array}{cc}
0&p\\
0&R
\end{array}
\right)$
almost does not depend on $p$. 

Here {\em ``almost does not depend''} means just {\em ``does not depend''} if $R\ne 0$ and {\em ``does not depend for $p \ne 0$ ''} if $R=0$. 
In the general case there is a similar
 non-degenerate condition.

If we consider $q,p,R$ as a blocks of proper dimensions, we arrive to the general case.
The parameterization of $A$ is reduced to the parameterization of $R$ with given Jordan structure. 
Symbolically the final formula can be written as
$$\omega_A=\tr dp\wedge dq + \omega_R,
$$
where $\omega_A$ is a symplectic form on the ``given'' orbit that contains $A$, 
$\omega_R $ is a symplectic form on the orbit that contains $R$, 
it has the strictly smaller dimension.

In the case $2\times 2$ the step ``{\em parameterization of $R$}'' is not visible, we should ``parameterize''  
the $0$-dimensional orbit of the $1\times 1$ matrix $R$.

 \vskip 2ex
 All the procedure uses the operation of solving linear equations systems only, consequently it is rational.
 
 The whole atlas for the orbit consists of the maps in question numerated by the permutations of the basic vectors. 
 The transition functions are rational, they are given by the formulae for the parameterization of the already parameterized matrix, but 
 conjugated by the matrices of the permutations of $N$ basic vectors.
 \vskip 3ex
 
 After all me make a good point.
As was noted at the very begining the
canonically conjugated families the constructed Archimedes-type coordinates
have special dual structure. The family that is an analog of the angle is generated
by the projections to the Grassmanians. 
It has a natural {\em global} structure as a set of affine coordinates on the projective manifold.

The analog of the projection on the diameter of the sphere in the original Archimedes scheme 
form the second family of the constructed coordinates. It has no evident global structure.
At the same time it has a remarkable {\em local} structure. 
The elements of the matrix on the orbit depend on these coordinates {\em linearly}.
 
\section{Definitions and notations}

It is well known that any coadjoint orbit of a semisimple Lie group is equipped with a standard symplectic 
(Lie-Poisson-Kirillov-Kostant) form.


In the partial case of $\operatorname{GL}(N, \CC)$ there are simplifications, at least in a terminology.
We can use the widely known language of elementary linear algebra in spite of the much less known language of 
Hamiltonian systems on Lie-algebras \cite{Lerman,Dusa}.

We use the Zariski topology, {\em open} set means {\em Zariski-open} set, {\em closed} set means {\em algebraically closed} set.
 
Let us treat elements $A\in \operatorname{gl}(N, \CC)$ as linear transformations $\mathcal A \in \operatorname{End} V$ of some 
complex linear $N$-dimensional space $V\simeq \CC^N$, equipped with a basis $\mathbf e_v=  e^1, \dots , e^N$:
$$A\in \operatorname{gl}(N, \CC) \leftrightarrow
\mathcal A: \mathbf e_v x\to \mathbf e_vAx, \ \ 
\mathbf e_v x, \mathbf e_v Ax\in V, \ \ x, Ax\in \CC^N.
$$
We identify the Lie algebra $\operatorname{gl}(N)$ and its dual $\operatorname{gl}^*(N)$ using the non-degenerate pairing (scalar product)
$\operatorname{gl}(N, \CC) \times \operatorname{gl}(N, \CC) \to \CC: <A,B>\to \tr AB$.
The Lie group $\operatorname{GL}(N, \CC) \ni g$ acts on $A\in \operatorname{gl}(N, \CC)$ by usual similarity transformations $A\to g^{-1}Ag$ induced by
changes of basis $\mathbf e_v \to\mathbf e_v g$, consequently an orbit of the coadjoint action can be identified 
with a manifold of all matrices similar to each other in this case.
Let us choose one element of the orbit, say $J$ that is the Jordan normal form of the matrices from the orbit,
and denote
$$\mathcal O(J)=\bigcup_{g\in \operatorname{GL}(N, \CC)}g Jg^{-1}.
$$
It is the subject of our investigations. 

The canonical symplectic Lie-Poisson-Kirillov-Kostant form $\omega_{{\mathcal O}(J)}$ on the orbit can be introduced 
by the equality
\begin{equation}\label{eq1}
\omega_{\scriptscriptstyle {\mathcal O}(J)}(\xi_1,\xi_2)=\tr J[g^{-1}\dot{g}_1,g^{-1}\dot{g}_2],
\end{equation}
where the vectors $\xi_1,\xi_2$ are tangent to the trajectories $A_i(t)=g_i(t)Jg_i^{-1}(t), \ i=1,2$ that intersect
each other at $t=0$:
$$g_1(0)=g_2(0)=g, \dot{g}_i=\left. \frac{d}{dt}\right|_{t=0}g_i(t),
\dot A_i=\left.\frac{d}{dt}\right|_{t=0}A_i(t).
$$
We will use the following version (see \cite{Hitchin})
of the previous formula
\begin{equation}\label{eq1b}
\omega_{\scriptscriptstyle {\mathcal O}(J)}(\xi_1,\xi_2)=\tr 
(\dot{g}_1g^{-1})\dot A_2
\end{equation}

The following observation (see \cite{DM,BD}) forms a basement of the construction: the canonical symplectic structure on an orbit
and the hierarchic structure (\ref{eq9}) which I present below are coordinated.
\vskip 3ex

Let $K\subset V$ be a subspace. Denote by $V/K$ the factor-space. It is a linear space of the dimension 
$\dim V-\dim K=\dim V/K$.
The linear structure is inherited from any $\dim V/K$-dimensional subspace of $V$ which is transverse to $K$.
We will denote by $\operatorname{Pr}^{\parallel K}$ the projection
\begin{equation}\label{eq2}
V\stackrel{Pr^{\parallel K}}{\longrightarrow} V/K.\end{equation}
Space $V$ has a structure of a trivial fiber bundle. Its fibers are subspaces parallel to $K$.

Let $\mathcal A$ be a linear transformation of $V$ and let its 
eigenspace corresponding to the eigenvalue $\lambda_0$ be $K:=\ker (\mathcal A-\lambda_0 id)$. Let
$$0 < \dim K < \dim V . 
$$
The submanifold 
of all  such $\mathcal A\in \operatorname{End} V$ 
will be denoted by $\operatorname{End} V|_{\lambda_0, K}$:
$$\mathcal A\in \operatorname{End} V|_{\lambda_0, K} \Longleftrightarrow \ker (\mathcal A-\lambda_0 id)=K.
$$

Let $\lambda_0=0$, $K=\ker \mathcal A$. The transformation $\mathcal A$ has the same value on all $X\in V$ from one equivalence class $V/K$ 
 that means that there is a linear transformation $\tilde{\tilde{\mathcal A}}\in \operatorname{Hom}(V/K, V)$ such that
$$\mathcal A=(\operatorname{Pr}^{\parallel K})^*\tilde{\tilde{\mathcal A}}.
$$
Let us denote by $((\operatorname{Pr}^{\parallel K})^{*})^{-1}$ the corresponding map 
$\operatorname{End} V|_{0, K}\to \operatorname{Hom} (V/K,V)$:
\begin{equation}\label{eq3}
\tilde{\tilde{\mathcal A}}=((\operatorname{Pr}^{\parallel K})^{*})^{-1}\mathcal A \Longleftrightarrow
\mathcal A = (\operatorname{Pr}^{\parallel K})^{*} \tilde{\tilde{\operatorname A}}
\end{equation}


The space $V$ has the structure of the fiber bundle 
$V\stackrel{Pr^{\parallel K}}{\longrightarrow} V/K$ consequently $((\operatorname{Pr}^{\parallel K})^{*})^{-1}$ 
can be projected back to $V/K$ by $\operatorname{Pr}^{\parallel K}$ that gives some
$\tilde{\mathcal A}\in \operatorname{End} V/K$:
$$\operatorname{Pr}^{\parallel K}\circ ((\operatorname{Pr}^{\parallel K})^{*})^{-1}: 
\operatorname{End} V|_{0,K} \rightarrow \operatorname{End} V/K.
$$
\begin{notation}
Let $\pi$ denote $\operatorname{Pr}^{\parallel K}\circ ((\operatorname{Pr}^{\parallel K})^{*})^{-1}$.
\end{notation}
To reconstruct the initial $\mathcal A$ from ${   \tilde{\mathcal A}}=
 \pi\mathcal A$
we need to know the position of the $\mathcal A$-image 
 on the assigned fiber of 
$V\stackrel{Pr^{\parallel K}}{\longrightarrow} V/K$.

Any subspace $M: M\oplus K=V$ 
sets the isomorphism $V/K\simeq M$ and defines the structure of the direct product on $\operatorname{End} V|_{0,K} $: 
\begin{equation}\label{eq4}
\operatorname{End} V|_{0,K} \stackrel{\sim}{\to} \pi(\operatorname{End} V|_{0,K})  \times \operatorname{Hom} (V/K, K).
\end{equation}

\section{Filtration of orbit}

The area of our exploration will be a modification of the Jordan structure by the action of the projection $\pi$.

By the Jordan structure $\mathcal J$ of a transformation $\mathcal A$ we mean the set of the eigenvalues of $\mathcal A$
and the information about the Jordan chains corresponding to each eigenvalue, namely the number of the chains and their lengths. 
By $J$ we denote a matrix (the normal Jordan form of $\mathcal A$) of the transformation $\mathcal A$ in some basis 
collecting from the vectors of the Jordan chains with the structure $\mathcal J$. We will specify the order of the vectors later.

An important property of the projection $\pi=\operatorname{Pr}^{\parallel K}\circ ((\operatorname{Pr}^{\parallel K})^{*})^{-1}$
on the first Cartesian factor of the target of (\ref{eq4}) we serve as the theorem.

\begin{theorem} \label{th1}
Let $\mathcal A$ be a linear transformation with non-trivial kernel $K: \ 0 < \dim K < \dim V $.

The Jordan structure ${   \tilde{\mathcal J}}$ 
of ${   \tilde{\mathcal A}}:=\pi\mathcal A$
is defined by the Jordan structure $\mathcal J$ of $\mathcal A$, namely
\begin{itemize}
\item the Jordan chains corresponding to the non-zero eigenvalues for $\mathcal J$ and for $\tilde{\mathcal J}$ coincide.

\item the Jordan chains corresponding to the zero eigenvalue of $\tilde{\mathcal J}$ are in one-to-one correspondence
with those chains of $\mathcal J$ that have non-unit length. The chains of $   \tilde{\mathcal J}$ are one unit shorter
than corresponding chains of $\mathcal J$.

\item the chains of the unit length (without generalized eigenvectors) form the kernel of the map (projection) of the 
set of Jordan chains of $\mathcal J$ to the set of Jordan chains of $\tilde{\mathcal J}$.
\end{itemize}

\end{theorem}

\underline{Proof}

First of all let us note that by the definition of $((\operatorname{Pr}^{\parallel K})^{*})^{-1}$
$$b=\mathcal A a \Rightarrow b=(((\operatorname{Pr}^{\parallel K})^{*})^{-1} \mathcal A) (\operatorname{Pr}^{\parallel K} a)
$$
consequently
\begin{equation}\label{eq1a}
{   \tilde{b}}:=\operatorname{Pr}^{\parallel K} b=(\pi\mathcal A)
(\operatorname{Pr}^{\parallel K} a)={   \tilde{\mathcal A}}{   \tilde{ a}}
\end{equation}
It implies that the cyclic law of the construction of Jordan chains takes place.

Consider any Jordan basis of $V$ for $\mathcal A$, where the first $\dim K$ vectors $e^1, \dots , e^{\dim K}$ form a basis of $K=\ker A$.

Consider the projection of the remaining subset of the basic vectors $e^{\dim K+1}, \dots, e^{\dim V} $.
It is a linear-independent set of the vectors, otherwise some linear combination
$\sum_{k>\dim K}\alpha_k e^k$ would be a vector from the kernel $K=\cup_{\alpha_k} \sum_{k\le\dim K}\alpha_k e^k$. 
It contradicts with the linear independence of $e^k$.

The number of vectors in the set $e^{\dim K+1}, \dots ,e^{\dim V} $ is equal to $\dim V -\dim K$ that is the dimension of $V/K$,
consequently the projection of the set
$$e^{\dim K +1}, \dots , e^{\dim V}
$$
forms a basis of $V/K$.
\qed

\begin{note}
The corresponding transformation of the Jordan structures  may be thought as {\em the projection } 
of the Jordan structures induced by the projection of a linear transformation along its kernel,
see the definition on the page~\pageref{project}.
\end{note}

\begin{corollary}\label{cor1}
All non-zero projections  by $\operatorname{Pr}^{\parallel K}$ 
of vectors forming {\em any} Jordan basis of $V$ for $\mathcal A$ form a Jordan basis of $V/K$ for 
$\pi\mathcal A$.

\end{corollary}



\begin{corollary}\label{cor2}
The projection of the set of the Jordan bases for $\mathcal A$ to the set of the Jordan bases for $\pi \mathcal A$ is surjective, 
namely for any Jordan basis ${   \tilde{\mathbf e}}_{J}$ of $V/K$ for 
$\pi\mathcal A$
there exists such a Jordan basis $\mathbf e_{J}$ of $V$ for $\mathcal A$ that the non-zero projections of its vectors
form ${   \tilde{\mathbf e}}_{J}$.

\end{corollary}
\underline{Proof}

By the definition of $\pi\mathcal A$ 
the statement that ${   \tilde{\mathbf e}}_{J}$ form a Jordan basis of $V/K$ for 
$\pi\mathcal A$  is equivalent to
the existence of the pre-images, i.e. it is equivalent to the existence of the set of vectors of $V$ connected by the Jordan cyclic law 
the projections of which are vectors of the set ${   \tilde{\mathbf e}}_{J}$.

It follows from the formula (\ref{eq1a}) that the cyclic law takes place. We can start from the pre-images 
of starting vectors of the Jordan chains of ${   \tilde{\mathbf e}}_{J}$ and iterate the transformation $\mathcal A$ in $V$.
The projection gives the iterations of $\pi\mathcal A$.
We get the set of vectors in $V$ that can be complemented to the basis by the eigenvectors of $\mathcal A$ without pre-images
(i.e. without the generalized eigenvectors), in other words the set can be complemented by the vectors of Jordan chains of the unit lengths.

Only one thing has to be proved. It is the linear independence of the last non-zero iterations of $\mathcal A$.
These iterations are already trivial for $\pi\mathcal A$,
because their inverse images belong to the $\ker \mathcal A$.

The statement follows from the uniqueness of the Jordan form. 
The desired dimension of the envelope of the last iterations is an invariant-defined value, it is the dimension of the
$\im \mathcal A\cap \ker \mathcal A$. In other words it is the difference between the 
dimension of $\ker A$ and the number of Jordan chains of the unit lengths. 

This number does not depend on the bases which we use for calculation, consequently it coincides with the number 
which we have for the basis constructed as the projection of any Jordan basis for $\mathcal A$ using the previous corollary.
\qed
\vskip 2ex

Denote by $\operatorname{End}_{\mathcal J}V$ a submanifold of the transformations with a fixed Jordan structure $\mathcal J$.
It has a structure of the fiber-bundle
\begin{equation}\label{eq5}
\operatorname{End}_{\mathcal J} V \stackrel{\gamma }{\to}
G(n, V), \ n:=\dim\ker J 
\end{equation}
over the Grassmanian. 

The fiber $\gamma^{-1}(K)$ over any $K\in G(n, V)$ is formed by all $\mathcal A$ from 
$\operatorname{End}_{\mathcal J} V\cap \operatorname{End} V|_{0,K}$.

It follows from (\ref{eq4}) that a fiber is a subset of $\operatorname{Hom(V/K,K)}$.
The  Jordan structure of $\mathcal A$ obviously will not be changed if we add {\em any} vectors 
from $K$ to the images of all vectors of a Jordan basis of $V$ for $\mathcal A$ if two restrictions are satisfied:
\begin{itemize}
\item the images of the vectors from $K$  keep zero values, 


\item the last vectors of chains  form a basis of $K$.
\end{itemize}
For the chains corresponding to the non-zero eigenvalues and for the chains corresponding to the zero eigenvalue 
but with the unit lengths these restrictions are trivial. 

For the chains of the lengths longer than one we have just one non-degeneracy restriction, namely {\em the images of the generalized 
eigenvectors of the first order (next to the last vectors of the chains corresponding to the zero eigenvalue)
must complete the set of vectors from the chains of the unit length to the basis of $K$.} The following theorem has been proved:

\begin{theorem}\label{th2}
A fiber $(\pi)^{-1}{   \tilde{\mathcal A}}$
of the bundle
\begin{equation}\label{eq6}
\operatorname{End}_{\mathcal J}V|_{0,K}
\stackrel{ \pi}{\longrightarrow} 
\operatorname{End}_{{   \tilde{\mathcal J}}}V/K
\end{equation}
is isomorphic to the open subset of $\operatorname{Hom(V/K,K)}$:
\begin{equation}\label{eq7}
\mathcal A \in (\pi)^{-1}{   \tilde{\mathcal A}}
\Leftrightarrow \rank \mathcal A|_{\scriptscriptstyle {   \tilde{K}}}=\dim {   \tilde{K}},
\end{equation}
where ${   \tilde{K}}:=(\operatorname{Pr}^{\parallel K})^{-1}\ker {   \tilde{\mathcal A}}$ is an inverse image
of the kernel of ${   \tilde{\mathcal A}}$ under the projection $\operatorname{Pr}^{\parallel K}:V\to V/K$:
$$x\in {   \tilde{K}} \Leftrightarrow \operatorname{Pr}^{\parallel K}x\in \ker {   \tilde{\mathcal A}}
$$
\end{theorem}
\qed

\begin{note}\label{nt2}
In the case $\ker \mathcal A \cap \im \mathcal A=0$, i.e. if $\mathcal A$ has no generalized eigenvectors for the zero eigenvalue
\begin{equation}\label{eq8}
\operatorname{End}_{\mathcal J} V|_{0,K}\simeq \operatorname{End}_{{   \tilde{\mathcal J}}} V/K \times 
\operatorname{Hom} (V/K,K).
\end{equation}
\end{note}
\qed


\vskip 2ex

We see that the manifold $\operatorname{End}|_{\mathcal J}$ has a structure of a fiber-bundle over the Grassmanian $G(n, V)$,
where the fiber is in its turn the fibration described by 
the previous theorem i.e. by the equalities (\ref{eq6}), (\ref{eq7}). In the simplest case of absence of
generalized eigenvectors it is given by the equality (\ref{eq8}).
\vskip 2ex

It is evident that for any eigenvalue $\lambda_1$ we can make the same construction with $\mathcal A - \lambda_1 id_v$, where by $id_v$
we denoted the identical transformation in $V$.
We get the similar representation, but from all the eigenvalues of all the chains the value $\lambda_1$ will be subtracted. 


Let us add $\lambda_1 id_{v/k}$ back to the transformations of $V/K$ 
in order that restores the initial set of eigenvalues. We introduce a special notation for such transformations
of Jordan structures. The transformation $\mathcal J \to {   \tilde{ \mathcal J}}$ from the Theorem~\ref{th1} is the partial case,
when $\lambda_1=0$.


\begin{definition}\label{project}
The operation of
{\em the projection of the Jordan structure $\mathcal J$ along the eigenspace, corresponding to the eigenvalue $\lambda_1$}
is a transformation of the Jordan structure $\mathcal J$ to the following Jordan structure denoted by $\mathcal J\setminus \{\lambda_1\}$: 
\begin{itemize}
\item all the Jordan chains corresponding to all $\lambda_i\ne \lambda_1$ are the same for $\mathcal J$ and  for 
$\mathcal J\setminus \{\lambda_1\}$

\item if all the Jordan chains of $\mathcal J$ corresponding to $\lambda_1$ have the lengths equal to one, 
$\mathcal J\setminus \{\lambda_1\}$ has no chains corresponding to the eigenvalue $\lambda_1$,
it consists of all the Jordan chains of $\mathcal J$ corresponding to $\lambda_i\ne \lambda_1$.
\item if $\mathcal J$ contains the Jordan chains corresponding to $\lambda_1$ of the lengths longer than one, 
$\mathcal J\setminus \{\lambda_1\}$
has chains corresponding to $\lambda_1$.
In this case the Jordan chains of $\mathcal J\setminus \{\lambda_1\}$ corresponding to the eigenvalue $\lambda_1$ are in one-to-one
correspondence with the Jordan chains of $\mathcal J$ corresponding to $\lambda_1$ with the lengths longer than one.
They are one unit shorter.
\end{itemize}
\end{definition} 

Let us denote 
$$\mathcal J\setminus \{\lambda_1\lambda_2 \dots \lambda_k\}:=
(\dots ((\mathcal J\setminus \{\lambda_1\})\setminus \{\lambda_2\}) \dots )\setminus \{\lambda_k\}$$
Note that in the case of the presence of generalized eigenvectors, the set $\{\lambda_1,\lambda_2, \dots, \lambda_k\}$ 
may contains the corresponding eigenvalue several times. It means we may 
{\em project along the eigenspace, corresponding to one eigenvalue} several times,
but no more times than the length of the longest chain corresponding to this eigenvalue is.

Consider a set $\lambda'_1,\lambda'_2, \dots, \lambda'_{M}$
of numbers collected from the set of eigenvalues of $\mathcal J$, where each eigenvalue $\lambda_k$ is written
such a number of times that is the length of the longest Jordan chain corresponding to it.
Consider any $\mathcal A\in \operatorname{End}_{\mathcal J} V$. It defines a point $K_1$ of the Grassmanian:
$$K_1=\ker ( \mathcal A -\lambda'_1 id_v)\in G(n_1,V),\ \  n_1:=\dim K_1 ,
$$ 
and the linear transformation of $V/K_1$:
$$\mathcal A_1:=\lambda'_1id_{v/k_1}+\operatorname{Pr}^{\parallel K_1}\circ ((\operatorname{Pr}^{\parallel K_1})^{*})^{-1}(\mathcal A
-\lambda'_1id_{v})\in \operatorname{End}_{\mathcal J\setminus\{\lambda'_1\}} V/K_1,
$$
where the sub-index near $id$ indicates the space where it is defined.

Let us consider the number $\lambda'_2$ and $V/K_1=:V_1$ where $\mathcal A_1$ acts. 
Due to the Theorem~\ref{th1}, $\lambda'_2$ is the 
eigenvalue of $\mathcal A_1$, so we can make the same procedure.
We get $V_2:=(V/K_1)/K_2$ and $\mathcal A_2$:
$$K_2=\ker ( \mathcal A_1 -\lambda'_2 id_{v_1})\in G(n_2,V_1),\ \  n_2:=\dim K_2 ,
$$ 
$$\mathcal A_2:=\lambda'_2id_{v_1/k_2}+\operatorname{Pr}^{\parallel K_2}\circ 
((\operatorname{Pr}^{\parallel K_2})^{*})^{-1}(\mathcal A_1 -\lambda'_2id_{v_1})\in 
\operatorname{End}_{\mathcal J\setminus\{\lambda'_1\lambda'_2\}} V_2,
$$
and so on, up to  the last $\mathcal J\setminus\{\lambda'_1, \dots \lambda'_{M-1}\}$ 
 for which the corresponding transformation is proportional to identical
$$\mathcal A_{M-1}=\lambda_{M}id.
$$
Denote a transformation of $A_{k-1}$ to $A_k$ by $\pi_{\{\lambda'_k\}}$ and consider a hierarchy

$$\operatorname{End}_{\mathcal J} V\stackrel{\gamma_1}{\longrightarrow} G(n_1,V)\ni K_1$$
$$
\gamma_1^{-1}(K_1)\stackrel{\pi_{\{\lambda'_1\}}}{\longrightarrow } 
\operatorname{End}_{\mathcal J\setminus\{\lambda'_1\}} V_1 \stackrel{\gamma_2}{\longrightarrow} G(n_2,V_1)\ni K_2$$

$$\dots \ \ \ \ \ \ \ \ \dots \ \ \ \ \ \ \ \ \dots $$

\begin{equation}\label{eq9}
\gamma_k^{-1}(K_k)\stackrel{\pi_{\{\lambda'_k\}}}{\longrightarrow } 
\operatorname{End}_{\mathcal J\setminus\{\lambda'_1\dots\lambda'_k\}} V_k 
\stackrel{\gamma_{k+1}}{\longrightarrow} G(n_{k+1},V_k)\ni K_{k+1}
\end{equation}

$$\dots \ \ \ \ \ \ \ \ \dots \ \ \ \ \ \ \ \ \dots $$

\begin{eqnarray*}
\gamma_{M-2}^{-1}(K_{M-2})\stackrel{\pi_{\{\lambda'_{M-2}\}}}{\longrightarrow } 
\operatorname{End}_{\mathcal J\setminus\{\lambda'_1\dots \lambda'_{M-2}\}} V_{M-2}\\
\stackrel{\gamma_{M-1}}{\longrightarrow} G(n_{M-1},V_{M-2})\ni K_{M-1}
 \end{eqnarray*}
$$\gamma_{M-1}^{-1}(K_{M-1})\stackrel{\pi_{\{\lambda'_{M-1}\}}}{\longrightarrow } 
\operatorname{End}_{\mathcal J\setminus\{\lambda'_1\dots \lambda'_{M-1}\}} V_{M-1},
 $$
 where 
 $$
\gamma_k:\operatorname{End}_{\mathcal J\setminus\{\lambda'_1\dots \lambda'_{k-1}\}}V_{k-1}\longrightarrow G(n_k, V_{k-1})
 $$ 
 maps any transformation to its eigenspace corresponding to $\lambda'_k$. Subspace $K_k$ is any $n_k$-dimensional subspace
 of $V_{k-1}$, $V_k:=V_{k-1}/K_k$. 
Transformation $\pi_{\{\lambda'_k\}}$ of $A_{k-1}$ to $A_k$ is defined by
\begin{eqnarray*}
\pi_{\{\lambda'_k\}}(\mathcal A_{k-1}):=\lambda'_kid+\operatorname{Pr}^{\parallel K_k}\circ 
((\operatorname{Pr}^{\parallel K_k})^{*})^{-1}(\mathcal A_{k-1} -\lambda'_kid)\\
=:\mathcal A_k\in 
\operatorname{End}_{\mathcal J\setminus\{\lambda'_1 \dots\lambda'_k\}} V_k,
\end{eqnarray*}
the transformation $((\operatorname{Pr}^{\parallel K_k})^{*})^{-1}$ is defined by (\ref{eq3}).
 \vskip 3ex

\section{Matrix representation}

 If a basis in $V$ is fixed linear transformations of $V$ get a matrix representation that identify 
 $\operatorname{End}_{\mathcal J} V$ and the manifold $\mathcal O(J)$ of all matrices similar to a given $J$.

Consider hierarchy (\ref{eq9}). The basis in $V$ does not induce neither matrix representations nor identifications 
$$\operatorname{End}_{\mathcal J\setminus \{\lambda'_1 \dots\lambda'_k\}} V_k \leftrightarrow 
\mathcal O({\mathcal J\setminus \{\lambda'_1 \dots\lambda'_k\}})$$  
on the levels of (\ref{eq9}) automatically,
because a projection of a basis is not a basis.

Consider any ordering of the vectors of the given basis $\mathbf e_v$.
Denote by $E_k$ the envelope of the last $m_k:=\dim \im \mathcal A_k$ vectors of $\mathbf e_v$, it is some coordinate subspace.
The sequence of projections along eigenspaces maps this set of $m_k$ vectors to $V_k$.
For each $\mathcal A$ from some open subset of $\mathcal O (J)$ this $\dim V_k$ vectors form a basis of $V_k$.

For each ordering of vectors of $\mathbf e_v$ this process sets natural isomorphisms between the abstract linear spaces $V_k$ 
and the coordinate subspaces $E_k$. The isomorphisms are defined for some open subset of the orbit $\mathcal O(J)$.

On the other hand a projection of a full set of vectors is a full set, consequently
for any $\mathcal A$ from the orbit  we can put vectors of ${\mathbf e}_v$ in such an order that the
bases of all $V_k$ will be formed by the images of the last several vectors of $\mathbf e_v$.

\begin{proposition}
The covering of the whole orbit $\mathcal O (J)$ by the open domains numerated by the permutations of vectors of $\mathbf e_v$
has been constructed.
\end{proposition}
\
Let us fix some ordering and identify $V_k$ with the corresponding subspaces $E_k$ of $V$. 

\begin{note}We will not distinguish $V_k$ and $E_k$ from now.
\end{note}



Filtration (\ref{eq9}) defines the sequence of the transformations $\mathcal A_k$
of the coordinate subspaces $E_k\simeq V_k$. In the given bases of $E_k$ the hierarchy 
has the transparent matrix representation:

$$A=\left(
\begin{array}{cc}
\operatorname I&0\\
Q_1&\operatorname I
\end{array}
\right)
\left(
\begin{array}{cc}
\lambda'_1\operatorname I&P_1\\
0 &A_1
\end{array}\right)
\left(
\begin{array}{cc}
\operatorname I&0\\
Q_1 &\operatorname I
\end{array}\right)^{-1}
$$
$$A_1=\left(
\begin{array}{cc}
\operatorname I&0\\
Q_2&\operatorname I
\end{array}
\right)
\left(
\begin{array}{cc}
\lambda'_2\operatorname I&P_2\\
0 &A_2
\end{array}\right)
\left(
\begin{array}{cc}
\operatorname I&0\\
Q_2 &\operatorname I
\end{array}\right)^{-1}
$$

$$\dots \ \ \ \ \ \ \ \ \dots \ \ \ \ \ \ \ \ \ \dots
$$

\begin{equation}\label{eq10}
A_{k-1}=\left(
\begin{array}{cc}
\operatorname I&0\\
Q_{k}&\operatorname I
\end{array}
\right)
\left(
\begin{array}{cc}
\lambda'_{k}\operatorname I&P_{k}\\
0 &A_k
\end{array}\right)
\left(
\begin{array}{cc}
\operatorname I&0\\
Q_{k} &\operatorname I
\end{array}\right)^{-1}
\end{equation}

$$\dots \ \ \ \ \ \ \ \ \dots \ \ \ \ \ \ \ \ \ \dots
$$

$$A_{M-2}=\left(
\begin{array}{cc}
\operatorname I&0\\
Q_{M-1}&\operatorname I
\end{array}
\right)
\left(
\begin{array}{cc}
\lambda'_{M-1}\operatorname I&P_{M-1}\\
0 &\lambda_{M} \operatorname I
\end{array}\right)
\left(
\begin{array}{cc}
\operatorname I&0\\
Q_{M-1} &\operatorname I
\end{array}\right)^{-1}
$$
\vskip 3ex

 
 Consider one flight of the hierarchy.
 To simplify the notations and in order to reduce the number of indexes let us put
 $k=0$, i.e. we consider the first flight. For the same reasons we put $\lambda_1=0$, the eigenspace
 $K_1$ becomes $\ker A=:K$, $\dim K=:n$, $\dim \im A=:m$.
 
 The coordinate subspace enveloping last $m$ basic vectors we denoted by $E_1=:E$,
 let us denote the envelope of the first $n$ vectors by $F$: $V=F\oplus E$.

 The open subset where this expansion of $V$ takes place consists of all the transformations
 $$ A\in \operatorname{End}_{ J} V: \ \ \ker A \cap E=0.
 $$
 
 Let us  denote such subset of the orbit by $(\!(\ \mathcal O(J)\ )\!)_E$. The corresponding subset of the Grassmanian
 we denote by $(\!(\ G(n,V) \ )\!)_E$: 
 $$K\in (\!( \ G(n,V)\ )\!)_E \Leftrightarrow K\cap E=0.$$
Let $(e^1 \dots e^n)$ be the set of the first vectors of $\mathbf e_v$. It forms the basis of $F$.

Let us denote the projection on the subspace $L_1$ along the subspace $L_2$ by   $\operatorname{Pr}^{\parallel L_2}_{L_1}$. 
It is defined for any transversal subspaces: $L_1\oplus L_2=V$. 

The projection along $E$ sets the isomorphism 
between $\ker A $ and $F$. It is defined for the transformations $A$ from the domain $(\!(\ \mathcal O(J)\ )\!)_E$.
Let us introduce the standard coordinates on the corresponding subset of the Grassmanian (see \cite{Harris}):
$$K\leftrightarrow \operatorname{Pr}^{\parallel F}_E\circ \operatorname{Pr}^{\parallel E}_K (e^1 \dots e^n).
$$
It is the isomorphism between the open subset of the Grassmanian in question and the set of all ${m\times n}$ matrices.

For our aims it is more natural to forget about the specification of the bases on $E$ and $F$ 
and use the decomposition $V=F\oplus E$ only.
It gives the following isomorphism $(\!( \ G(n,V)\ )\!)_E\stackrel{\sim}{\to}\operatorname{Hom}(F,E)$:
\begin{equation}\label{eq11}
K\leftrightarrow \operatorname{Pr}^{\parallel F}_E\circ \operatorname{Pr}^{\parallel E}_K|_F 
\end{equation}

Consider the fibration that is given by Theorem~\ref{th2}. Its fiber is the open subset
of $\operatorname{Hom}(E,K)$. We know that the projection parallel to $E$ sets the
bijection between $F$ and $K=\ker A$,  on $(\!(\ \mathcal O(J)\ )\!)_E$, 
consequently we may replace $\operatorname{Hom}(E,K)$
on $\operatorname{Hom}(E,F)$. We should just compose the projection $\operatorname{Pr}^{\parallel E}_F$ with
each element of $\operatorname{Hom}(E,K)$.

We formulate the modified version of the Theorem~\ref{th2} as the proposition.

\begin{proposition}\label{prop1}
The open set $(\!(\ \mathcal O(J)\ )\!)_E$ has a structure of a fiber-bundle
$$(\!(\ \mathcal O(J)\ )\!)_E\to \operatorname{Hom}(F,E)\times \mathcal O(J\setminus\{0\}),
$$
where a fiber is the open subset of $\operatorname{Hom}(E,F)$.
\end{proposition}
\vskip 3ex

Theorem~\ref{th2} states that some algebraically closed set of $\operatorname{Hom}(E,F)$
does not belong to the fiber. It corresponds to the transformations which Jordan chains are shorter than
necessary and their last units do  not form a basis of the eigenspace.

This effect comes into particular prominence in $2\times 2 $ case. Let us parameterize the Zariski-open part
$(\!(\ \mathcal O(J)\ )\!)_E$ of the orbit\footnote{
In this case the complement of $(\!(\ \mathcal O(J)\ )\!)_E$ to $(\!(\ \mathcal O(J)\ )\!)$ is formed
by the lower triangular matrices  $\left( \begin{array}{cc}0&0\\ x &0\end{array}\right)$. The subspace $E$ is
the coordinate subspace spanned the second basic vector $e^2$} of 
$J=\left( \begin{array}{cc}0&1\\ 0&0\end{array}\right)$ by the functions $p,q$:
$$\left( \begin{array}{cc}1&0\\ q&1\end{array}\right)\left( \begin{array}{cc}0&p\\ 0&0\end{array}\right)
\left( \begin{array}{cc}1&0\\ q&1\end{array}\right)^{-1}
=\left( \begin{array}{cc}-pq&p\\ -pq^2 &pq\end{array}\right).
$$
We must exclude $p=0$ because $\left( \begin{array}{cc}0&0\\ 0&0\end{array}\right)$ does not belong to the orbit.

Nevertheless in several important applications it is not naturally to exclude the divisor $p=0, q\in \CC$ from the chart.
For example for the problems of isomonodromic deformations the corresponding Hamiltonian flows have no any special
behavior on the line $p=0, q=\CC$, so the natural way to make a theory consistent is {\em to add this divisor} i.e.
to expand the orbit $\mathcal O (J)$. 
It means that we introduce a new symplectic manifold $\mathcal O'(J)$ with almost
the same fiber-bundle structure as $\mathcal O (J)$ has, but the fibers are just $\operatorname{Hom}(E,F)$.

In other words to get $\mathcal O' (J)$ we glue some algebraically closed set to $\mathcal O (J)$ in such a way
that $\mathcal O' (J)$ will be the symplectic manifold too and there will be the symplectic map 
$\mathcal O (J)\hookrightarrow\mathcal O' (J)$. 

\begin{note}\label{nt3}
We assume that our orbits are already enlarged. We mark this enlargering by the prime. It means
that we are investigating the symplectic manifolds $\mathcal O' (J)$ that are equipped with the symplectic 
mappings $\mathcal O' (J)\hookleftarrow \mathcal O (J)$ .

The complement to the image of $\mathcal O (J)$ in $\mathcal O' (J)$ is the algebraically closed 
set isomorphic to the set of the transformations which Jordan chains (their last units)
do not form a basis of the corresponding eigenspace.
\end{note}

From now we are in the conditions of the Note~\ref{nt2}. We formulate the version of the Theorem~\ref{th2} for the
enlarged orbits $\mathcal O'$ as the next proposition.

Coordinate subspace  $E_k$ is identified with $V_k$. 
In accordance to the identification let us change the notation of
 $((\operatorname{Pr}^{\parallel K})^{*})^{-1}$ to $((\operatorname{Pr}^{\parallel K_k}_{E_k})^{*})^{-1}$. 
We keep the previous notation $\pi_{\{\lambda'_k\}}$ for the transformation $A_{k-1}\to A_k$.

\begin{proposition}\label{prop2}
The open set $(\!(\ \mathcal O'(J\setminus\{\lambda'_1\dots \lambda'_{k-1}\})\ )\!)_{E_k}$ 
is isomorphic to the direct product
$$(\!(\ \mathcal O'(J\setminus\{\lambda'_1\dots \lambda'_{k-1}\})\ )\!)_{E_k}
\stackrel{\sim}{\to} \operatorname{Hom}(E_k,F_k)\times\operatorname{Hom}(F_k,E_k)\times \mathcal O(J\setminus\{\lambda'_1\dots \lambda'_k\})
$$
A point $A_{k-1}\in (\!(\ \mathcal O(J\setminus\{\lambda'_1\dots \lambda'_{k-1}\})\ )\!)_{E_k}$ has the 
following projections on the Cartesian factors
$$A_{k-1}\to (\mathcal P_k,\mathcal Q_k, A_k):
$$
$$\mathcal P_k=\operatorname{Pr}^{\parallel E_k}_{F_k}\circ 
((\operatorname{Pr}^{\parallel K_k}_{E_k})^{*})^{-1}(A_{k-1} -\lambda'_k\operatorname{I}),$$
$$\mathcal Q_k=\operatorname{Pr}^{\parallel F_k}_{E_k}\circ \operatorname{Pr}^{\parallel E_k}_{K_k}|_{F_k},$$
$$A_k=\lambda'_k\operatorname{I}+\operatorname{Pr}^{\parallel K_k}_{E_k}\circ 
((\operatorname{Pr}^{\parallel K_k}_{E_k})^{*})^{-1}(A_{k-1} -\lambda'_k\operatorname{I})=:\pi_{\{\lambda'_k\}}A_{k-1}.
$$


\end{proposition}
\begin{note}
The formulae for the projections have sense for the points of $\mathcal O$ only.
The complement of the 
$\mathcal O(J\setminus\{\lambda'_1\dots \lambda'_{k-1}\})$ with respect to the
$\mathcal O'(J\setminus\{\lambda'_1\dots \lambda'_{k-1}\})$  does not immersed into $\operatorname{gl}(N_{k-1})\simeq \operatorname{End} V_{k-1}$.
\end{note}

\begin{note}
The proposition is equivalent to (\ref{eq10}) where the restrictions on the matrices $P_k$ are removed.
\end{note}

\subsection{Auxiliary symplectic space $\mathcal E(F\oplus E)$}

Let $V=F\oplus E$  be a splitting of $n+m$-dimensional space $V$ into 
a direct sum of two subspaces $F$, and $E$, \ $ \dim F=n,\dim E=m$.

Consider $\operatorname{End} V$ as a linear $(n+m)^2$-dimensional space.
Let us define a skew-symmetrical scalar product $\omega_{\mathcal E}: \operatorname{End} V\times \operatorname{End} V\to \CC$:
\begin{equation}\label{eq12}
\omega_{\mathcal E}(\mathcal B_1, \mathcal B_2)=
\tr \operatorname{Pr}^{\parallel F}_{E}\mathcal B_2\circ \operatorname{Pr}^{\parallel E}_{F} \mathcal B_1|_E
- \tr \operatorname{Pr}^{\parallel E}_{F}\mathcal B_2\circ \operatorname{Pr}^{\parallel F}_{E} \mathcal B_1|_F.
\end{equation}
It is obviously degenerated.
Let us introduce a basis $(\mathbf{fe})$, where $\mathbf{f}$ is a basis of $F$ and $\mathbf{e}$ is a basis of $E$.
Let $\mathcal B_i$ have the matrix $\left( \begin{array}{cc}b^i_{11}&b^i_{12}\\ b^i_{21}&b^i_{22}\end{array}\right)$ 
in this basis. By $b^i_{jk}$ we denoted the corresponding blocks.

In this coordinates 
\begin{equation}\label{eq13}
\omega_{\mathcal E}(\mathcal B_1, \mathcal B_2)=\tr b^2_{21}b^1_{12} - \tr b^2_{12}b^1_{21}.
\end{equation}

Let us introduce the $2nm$-dimensional subspace $\mathcal E(F\oplus E)\subset \operatorname{End} V$:
$$\mathcal B\in \mathcal E(F\oplus E)\Leftrightarrow \operatorname{Pr}^{\parallel F}_{E} \mathcal B|_E=0
=\operatorname{Pr}^{\parallel E}_{F} \mathcal B|_F.
$$
It consists of the matrices which are 
off-diagonal in the basis $(\mathbf{fe})$, i.e. $b_{11}=0=b_{22}$. We keep the previous notation
$\omega_{\mathcal E}$ for the restriction of $\omega_{\mathcal E}$ on $\mathcal E(F\oplus E)$
.

\begin{proposition}\label{prop3}
Space $(\mathcal E(F\oplus E), \omega_{\mathcal E}) $ is $2nm$-dimensional symplectic space.
Canonical basis is formed by the set of couples $P_{ij}, Q_{ji}, \ 1\le i\le n, n+1\le j\le n+m$, where
$P_{ij}, Q_{ji}\in \mathcal E(F\oplus E)$ are the transformations with the following matrices
$$(P_{ij})_{st}=\delta_{si}\delta_{tj}, \ \ \  (Q_{ji})_{st}=\delta_{sj}\delta_{ti}
$$
in the basis $(\mathbf{fe})$ where the first $n$ vectors form basis $F$ and the last $m$ vectors form basis $E$.
\end{proposition}
\underline{Proof}

The proof follows from the formula (\ref{eq13}) that shows that $P_{ij}, Q_{ji}$ is really the Darboux basis
for $(\mathcal E(F\oplus E), \omega_{\mathcal E}) $.
\qed

\begin{proposition}\label{prop4}
There is a natural isomorphism between the manifold $\operatorname{Hom}(E,F)\times\operatorname{Hom}(F,E)$
and the space $\mathcal E(F\oplus E)$.
\end{proposition}

To construct the point of $\operatorname{End} V\supset \mathcal E(F\oplus E)$ 
 it is sufficient to assign its action on each of summands of $F\oplus E$. 
Let a couple $\mathcal P,\mathcal Q$ be a point of $\operatorname{Hom}(E,F)\times\operatorname{Hom}(F,E)$.
We define the  transformation $\mathcal B\in \operatorname{End} V$  corresponding to the couple as the map
which transforms the vectors from $E,F\subset V$ (the natural embedding) as it is assigned by $\mathcal P$ and $\mathcal Q$.

Consider the opposite direction. Any $\mathcal B\in \operatorname{End} V$ can be decomposed on 
$\mathcal P:=\operatorname{Pr}^{\parallel E}_{F}\mathcal B|_E$ and
$\mathcal Q:=\operatorname{Pr}^{\parallel F}_{E}\mathcal B|_F$.
For the transformations from $\mathcal E(F\oplus E)\subset
\operatorname{End} V$ these $\mathcal P$ and $\mathcal Q$ define $\mathcal B$ uniquely. 
It is evidently an isomorphism.

\qed
\vskip 3ex
\subsection{Main theorem}

We can see that out of some algebraically closed subset the symplectic manifold $\mathcal O(J)$
is isomorphic to the Cartesian product of the  {\em linear} symplectic space 
$\mathcal E(F\oplus E)$ and the symplectic manifold $\mathcal O (J\setminus \{\lambda'\})$
of the smaller than $\mathcal O(J)$ dimension. Here $\lambda'$ is some eigenvalue of $J$ and
$$F\simeq \ker (J-\lambda' \operatorname{I}) , \ \ \ E\simeq \im (J-\lambda' \operatorname{I}).
$$

Let us denote the projection on the Cartesian factor $\mathcal E(F\oplus E)$ by 
$\pi_{\scriptscriptstyle \mathcal E(F\oplus E)}$
and the projection on $\mathcal O(J\setminus \{\lambda'\})$ by $\pi_{\scriptscriptstyle \{\lambda'\}}$.

We constructed the isomorphism between two symplectic spaces
equipped with their own forms  $\omega_{\scriptscriptstyle \mathcal O(J)}$ and 
$\omega_{\scriptscriptstyle \mathcal E(F\oplus E)}+ \omega_{\mathcal O(J\setminus\{\lambda'\})}$.


Let us introduce the main theorem now. 

\begin{theorem}\label{main}
The isomorphism
\begin{equation}\label{eq15}
(\!( \ \mathcal O'(J)\ )\!)_{E}
\stackrel{\sim}{\to} \mathcal E(F\oplus E) \times \mathcal O(J\setminus\{\lambda'\})
\end{equation}
is birational and symplectic:
\begin{equation}\label{eq17}
\omega_{\scriptscriptstyle \mathcal O(J)}=
\pi^*_{\scriptscriptstyle \mathcal E(F\oplus E)}\omega_{\scriptscriptstyle \mathcal E(F\oplus E)}+ 
\pi^*_{\scriptscriptstyle \{\lambda'\}}\omega_{\scriptscriptstyle\mathcal O(J\setminus\{\lambda'\})}.
\end{equation}
\end{theorem}

The proof will be based on the following lemma.
\vskip 2ex

Let the given basis $\mathbf e_v$ be divided in two parts $\mathbf e=(\mathbf f{   \tilde{\mathbf e}})$
in accordance with the dimensions of the kernel and the image of $A-\lambda' \operatorname{I}$.
Let $E$ be the envelope of $\tilde{\mathbf e}$, $A\in (\!( \ \mathcal O(J)\ )\!)_{E}$,  
${   \tilde{A}}:=\pi_{\scriptscriptstyle \{\lambda'\}}A\in\mathcal O(J\setminus\{\lambda'\})$, 
$\dim\ker (A-\lambda' \operatorname{I})=n$,  $\dim\im (A-\lambda' \operatorname{I})=m$.

\begin{lemma}\label{lem1}

For any $g\in \operatorname{GL}(m, \CC)$ that transforms 
the fixed basis ${   \tilde{\mathbf e}}$ of $E$ to any Jordan basis ${   \tilde{\mathbf e}}_J$ of 
$E$ for ${   \tilde{A}}$,
 there exist 
 \begin{itemize}
 \item the the set of vectors $\kappa$ that form the basis of $\ker (A-\lambda' \operatorname{I})$, 
\item the matrix $\hat P\in\CC^{n\times m}$ 
\end{itemize}
 such that $\mathbf e_J$:
\begin{equation}\label{eq18}
\mathbf e_J=
(\mathbf{\kappa }{   \tilde{\mathbf e}})
\left(
\begin{array}{cc}\operatorname{I} &\hat P\\0&g \end{array}
\right)=
(\mathbf{ \kappa} {   \tilde{\mathbf e}}_J)
\left(\begin{array}{cc}\operatorname{I} &\hat P\\0&\operatorname{I} \end{array}
\right)
\end{equation}
 form a Jordan basis of $V$ for $A$. 
\end{lemma}

\underline{Proof of the lemma}

For the simplification of the notations let us put $\lambda'=0$. Consider $\mathcal O(J)$, where $J$
is the Jordan normal form of the matrices from the orbit. 
Let us order the vectors of the Jordan basis for $J$ in 
such a way that the first set $\kappa$ of the vectors of the basis of $V$ 
forms the basis of the root-space of $J$:
$$J=\left( \begin{array}{cc} 0&J_P\\ 0& {   \tilde{J}} \end{array}\right).
$$
Consider such a part of the orbit where $\kappa$ is completing 
some fixed linear independent set ${   \tilde{\mathbf e}}$
to the basis of $V$.

In the basis $(\kappa, {   \tilde{\mathbf e}})$ any $A$ from the orbit has the form
 $$\left(  \begin{array}{cc} 0&P\\ 0& {   \tilde{A}} \end{array}\right).
 $$
 The statement of the lemma is equivalent to the following:
 
 {\em if $\left(  \begin{array}{cc} 0&P\\ 0& {   \tilde{A}} \end{array}\right)$ is similar to the
 $\left(  \begin{array}{cc} 0&J_P\\ 0& {   \tilde{J}} \end{array}\right)$, and if the zero columns form the basis
 of the root-spaces of the matrices,
 
 then for the given $g$: $g^{-1}{   \tilde{A}}g={   \tilde{J}}$ there exist such $\hat g\in\operatorname{GL}(n, \CC)$, 
 and such $\hat P$ that
 \begin{equation}\label{lemma}
 \left(\begin{array}{cc} 0&P\\ 0& {   \tilde{A}}\end{array}\right)=
 \left(
\begin{array}{cc}\hat g&\hat P\\0&g \end{array}\right)
 \left(  \begin{array}{cc} 0&J_P\\ 0& {   \tilde{J}}\end{array}\right)
 \left( \begin{array}{cc}\hat g&\hat P\\0&g \end{array}\right)^{-1}.
 \end{equation}
 }
 
 It is equivalent to the solvability of the equation on $\hat g$ and $\hat P$
 $$P=(\hat g J_P+\hat P {   \tilde{J}})g^{-1},
 $$
 where $P,g,J_P,{   \tilde{J}}$ are given.

The equation is solvable for any $P$, because the number of the linear independent rows 
in $(n+m)\times m$ matrix $\left(  \begin{array}{c}J_P\\ {   \tilde{J}} \end{array}\right)$
coincides with the number of the linear independent columns that is $m$. It is just the dimension of the
the space $\CC^{m}$ of  the rows of $P$.

To prove that the matrix $\hat g$ can be chosen non-degenerated 
let us rewrite the equation:

%

\begin{equation} \label{eq23}
Pg=\hat g J_P+\hat P {   \tilde{ J}}.
\end{equation}
 
The matrices $J_P$ and $\tilde J$ are the blocks of the Jordan matrix 
$J=\left( \begin{array}{cc} 0&J_P\\ 0& {   \tilde{J}} \end{array}\right)$, each column of 
$\left( \begin{array}{c} J_P\\ \tilde{J} \end{array}\right)$ contains exactly one unit.
It implies that the root space of ${   \tilde{ J}}$ and the root space of $J_P$ form a basis of columns $\CC^m$.
Consider the zero columns of the Jordan matrix ${   \tilde{ J}}$. The set of the corresponding columns of $Pg$
has full dimension otherwise there will be a linear relation between the columns of 
 $\left(  \begin{array}{c}P\\ {   \tilde{A}} \end{array}\right)$:
 $$P-\hat P {   \tilde{J}}g^{-1}=P-\hat Pg^{-1} {   \tilde{ A}}=\hat gJ_Pg^{-1}. $$
 Consider (\ref{eq23}).
 From the linear independence of the columns of $Pg$ in question it follows  
 that on the places of the zero columns of ${   \tilde{J}}$ there are linear
 independent columns of $Pg$ that implies the linear independence of the corresponding columns of $\hat g$.
 Matrix $\hat gJ_P$ does not depend on the other columns of $\hat g$ because the 
 corresponding columns of $J_P$ vanishes, consequently
 the set of the linear independent columns of $\hat g$  can be completed in an arbitrary way, 
 we choose $\det \hat g\ne 0$.

\underline{The lemma has been proved}
\vskip 2ex

\begin{note}
Lemma itself follows from the Corollary~\ref{cor2} directly, but for the future considerations we need the information
about the introduced matrices. 
\end{note}
 Let us proof the theorem.
 \vskip 3ex
 
 Consider any point $A\in (\!( \ \mathcal O(J)\ )\!)_{E}$, and the level-sets of the map (\ref{eq15}):
 $(\cup A)|_{\mathcal E=\const}$ and $(\cup A)|_{{{\mathcal O}}=\const}$ passing this point.
 The map (\ref{eq15}) is the isomorphism, consequently
 $$\operatorname{T}_A\mathcal O(J)=
 \operatorname{T}_{A}(\cup A)|_{\mathcal E=\const}
 \oplus \operatorname{T}_{A}(\cup A)|_{{{\mathcal O}}=\const}.
 $$

Let $\partial_{\mathcal E}$ and $\partial_{\mathcal O}$ be any vectors from the corresponding subspaces:
$$\partial_{\mathcal E}\in \operatorname{T}_{A}(\cup A)|_{\mathcal O=\const}, \ \
\partial_{\mathcal O}\in\operatorname{T}_{A}(\cup A)|_{{{\mathcal E}}=\const}.
$$

They are tangents to the lines
$$A_{\mathcal O}(t)=
\left(\begin{array}{cc} \operatorname{I}&0\\ Q(t)&\operatorname{I}\end{array}\right)
  \left(\begin{array}{cc}0&P(t)\\0&{   \tilde{A}} \end{array}\right)
  \left(\begin{array}{cc} \operatorname{I}&0\\ Q(t)&\operatorname{I}\end{array}\right)^{-1}
+\lambda'\operatorname{I}.
$$
and
$$A_{\mathcal E}(t)=
\left(\begin{array}{cc} \operatorname{I}&0\\ Q&\operatorname{I}\end{array}\right)
  \left(\begin{array}{cc}0&P\\0&{   \tilde{A}}(t) \end{array}\right)
  \left(\begin{array}{cc} \operatorname{I}&0\\ Q&\operatorname{I}\end{array}\right)^{-1}
+\lambda'\operatorname{I}
$$
that belong to the corresponding level sets.

It follows from the lemma that on the level set $(\cup A)|_{\mathcal O=\const}$ any curve $A_{\mathcal O}(t)$
can be parameterized in the following way:
$$\left(\begin{array}{cc} \operatorname{I}&0\\ Q(t)&\operatorname{I}\end{array}\right)
 \left(\begin{array}{cc}\hat g(t)&\hat P(t)\\0& g\end{array}\right)
  \left(\begin{array}{cc}0&J_P\\0&{   \tilde{J}} \end{array}\right)
 \left(\dots 
 \right)^{-1}
\left(  \dots
\right)^{-1}+\lambda'\operatorname{I},
$$
consequently
$$\left.\frac{d}{dt}\right|_AA_{\mathcal O}(t)=
\left[\left(\begin{array}{cc} \operatorname{I}&0\\ Q&\operatorname{I}\end{array}\right)
\left(\begin{array}{cc} *&*\\ *&0\end{array}\right)
\left(\begin{array}{cc} \operatorname{I}&0\\ Q&\operatorname{I}\end{array}\right)^{-1},
A\right]
$$
$$\left.\frac{d}{dt}\right|_AA_{\mathcal E}(t)=
\left(\begin{array}{cc} \operatorname{I}&0\\ Q&\operatorname{I}\end{array}\right)
  \left(\begin{array}{cc}0&0\\0&* \end{array}\right)
  \left(\begin{array}{cc} \operatorname{I}&0\\ Q&\operatorname{I}\end{array}\right)^{-1},
$$
where we denote the terms the values of which are unessential by stars.

The application of formula (\ref{eq1b}) gives the desired
\begin{equation}\label{eq19}
\omega_{\scriptscriptstyle {\mathcal O}(J)}(\partial_{\mathcal E},\partial_{\mathcal O})=0
\end{equation}
Let $\partial^i_{\mathcal Q}\in \operatorname{T}_{A}(\cup A)|_{\mathcal O=\const}, \ i=1,2$ be two vectors tangent
to the level-set of function $P$ i.e. they are tangents to the lines
$$A^i_{\mathcal{PO}}(t)=
\left(\begin{array}{cc} \operatorname{I}&0\\ Q_i(t)&\operatorname{I}\end{array}\right)
  \left(\begin{array}{cc}0&P\\0&{   \tilde{A}} \end{array}\right)
  \left(\begin{array}{cc} \operatorname{I}&0\\ Q_i(t)&\operatorname{I}\end{array}\right)^{-1}
+\lambda'\operatorname{I}.
$$
The calculation gives:
$$ \omega_{\scriptscriptstyle {\mathcal O}(J)}(\partial^1_{\mathcal Q},\partial^2_{\mathcal Q})
=\tr
\left(
\begin{array}{rr}
0&0\\
{\dot{Q}}_1&0
\end{array}\right)\left[\left(
\begin{array}{rr}
0&0\\
{\dot{Q}}_2&0
\end{array}\right), A\right]=0.
$$

Let $\partial^i_{\mathcal P}\in \operatorname{T}_{A}(\cup A)|_{\mathcal O=\const}, \ i=1,2$ be two vectors tangent
to the level-set of function $Q$ i.e. they are tangents to the lines
$$A^i_{\mathcal{QO}}(t)=
\left(\begin{array}{cc} \operatorname{I}&0\\ Q&\operatorname{I}\end{array}\right)
  \left(\begin{array}{cc}0&P_i(t)\\0&{   \tilde{A}} \end{array}\right)
  \left(\begin{array}{cc} \operatorname{I}&0\\ Q&\operatorname{I}\end{array}\right)^{-1}
+\lambda'\operatorname{I}.
$$
We can set $Q=0$ because trace does not depend on the conjugation of all the factors by one matrix.
The calculation of the tangent vectors gives
 $$\left.\frac{d}{dt}\right|_A
 \left(\begin{array}{cc} 0&P_i(t)\\ 0& {   \tilde{A}}\end{array}\right)=
 \left(\begin{array}{cc} 0&\dot P_i\\ 0& 0 \end{array}\right)$$
 $$=
 \left.\frac{d}{dt}\right|_A
 \left(\begin{array}{cc}\hat g_i(t)&\hat P_i(t)\\0&g \end{array}\right)
 \left(  \begin{array}{cc} 0&J_P\\ 0& {   \tilde{J}}\end{array}\right)
 \left( \begin{array}{cc}\hat g_i(t)&\hat P_i(t)\\0&g \end{array}\right)^{-1}$$
 $$=
 \left[\left(\begin{array}{cc}*&*\\0&0 \end{array}\right),
 \left(\begin{array}{cc} 0&P\\ 0& {   \tilde{A}} \end{array}\right) \right]
 $$
Consequently
$$ \omega_{\scriptscriptstyle {\mathcal O}(J)}(\partial^1_{\mathcal P},\partial^2_{\mathcal P})
=\tr
\left(
\begin{array}{rr}
0&\dot P_1\\
0&0
\end{array}\right)\left(
\begin{array}{rr}
*&*\\
0&0
\end{array}\right)=0.
$$

Simple calculation gives 
$\omega_{\scriptscriptstyle {\mathcal O}(J)}(\partial_{\mathcal P},\partial_{\mathcal Q})=\tr PQ$,
that means 
\begin{equation}\label{eq21}\left.\omega_{\scriptscriptstyle \mathcal O(J)}\right|_{\mathcal O=\const}=
\pi^*_{\scriptscriptstyle \mathcal E(F\oplus E)}\omega_{\scriptscriptstyle \mathcal E(F\oplus E)}.
\end{equation}

Let us consider two tangents $\partial^1, \partial^2$ to the lines $A_i(t)$ on the level-set $(\cup A)|_{\mathcal E=\const}$. 
For the previous reasons without the loss of generality we put $Q=0$,
$$\left.\frac{d}{dt}\right|_A
 \left(\begin{array}{cc} 0&P\\ 0& {   \tilde{A}}_i(t)\end{array}\right)=
 \left(\begin{array}{cc} 0&0 \\ 0& \dot{{   \tilde{A}}}_i \end{array}\right).$$
 From the representation (\ref{lemma}) we get
 $$\left.\frac{d}{dt}\right|_A \left(\begin{array}{cc} 0&P\\ 0& {   \tilde{A}}_i(t)\end{array}\right)=
 \left[
 \left(\begin{array}{cc} *&*\\ 0& g^{-1}\dot g_i\end{array}\right)
 ,\left(\begin{array}{cc} 0&P\\ 0& {   \tilde{A}}\end{array}\right)
 \right].$$
The application of the formula (\ref{eq1b}) gives
 $$\omega_{\scriptscriptstyle \mathcal O(J)}(\partial^1,\partial^2)=
 \tr g^{-1}\dot g_1 \dot{{   \tilde{A}}}_2,
$$
that is the value of $\omega_{\scriptscriptstyle \mathcal O(\tilde J)}$ on the projections of the vectors $\partial^1,\partial^2$.

The equality (\ref{eq17}) follows from (\ref{eq19}), (\ref{eq21}) and the last one.
\vskip 2ex

 Let us prove the birationality of the isomorphism (\ref{eq15}). 
To find the images of the projections 
$ \pi_{\scriptscriptstyle \mathcal E(F\oplus E)}$ and 
$\pi_{\scriptscriptstyle \{\lambda'\}}$
 we have to find the eigenvectors corresponding to the given eigenvalue and project along the subspaces.
 The inverse operation is the multiplication of the matrices with the given blocks in formulae (\ref{eq10}).
All these operations are rational.
 
 \qed

\vskip 3ex

Let us present the final formulae for the map
$$\mathcal E(F_1\oplus E_1)\times \mathcal E(F_2\oplus E_2)\times \dots \times \mathcal E(F_M\oplus E_M) 
\to \mathcal O'(J).
$$

Let $\lambda'_1, \dots , \lambda'_M$ be a sequence of the eigenvalues of some Jordan matrix
matrix $J$. Let each eigenvalue $\lambda'$ be written in the sequence such a number of times $r_{\lambda'}$ as the length
of the longest Jordan chain corresponding to this eigenvalue is:
$$\dim\ker (J-\lambda' I)^{r_{\lambda'}-1}<\dim\ker (J-\lambda' I)^{r_{\lambda'}}=\dim\ker (J-\lambda' I)^{r_{\lambda'}+1}.
$$
Let us denote by $n_k$ 
{\em a number of the Jordan chains  
that are not shorter than the number of eigenvalues equal to this $\lambda'_k$ in the subsequence $\lambda'_1, \dots, \lambda'_k$ }.

\begin{proposition}
The  full information about the Jordan structure of $J$ is contained in the set of couples $(\lambda'_k,n_k), \ k=1,\dots, M$.
\end{proposition}

\underline{Proof}

If the eigenspace corresponding to $\lambda'$ does not contain generalized eigenvectors we have $r_{\lambda'}=1$. 
In this case there is 
only one $\lambda'_k=\lambda'$ in the sequence and the number $n_k$ is the dimension of the eigenspace.

Let the set of Jordan chains corresponding to $\lambda'$  consists of
$m_1$ chains of the length $1$, $m_2$ chains of the length $2$, \dots , $m_{r_{\lambda'}}\ne 0$ chains of the length $r_{\lambda'}$.
In this case the set of the numbers $n_k$ corresponding to these eigenvalue is the non-increasing sequence of $r_{\lambda'}$  integers 
$m_i+m_{i+1}+\dots+m_{r_{\lambda'}}, \ i=1,\dots r_{\lambda'}$. The smallest $n_k$ is the number of the longest Jordan chains. Their lengths 
$r_{\lambda'}$ are equal to the number of repetitions of $\lambda'$ in the sequence. 
To reconstruct other $m_i$'s we should take the differences between neighbour $n_k$'s.
\qed

Denote by $Q$ the following lower-triangular block-matrix
$$Q:=\left(
\begin{array}{ccccccc}
I_{n_1} &0 &0&\dots &0&0&0\\
q_2^1&I_{n_2} &0 &\dots &0&0&0\\
q_3^1&q_3^2&I_{n_3} &\dots &0&0&0\\
\dots&\dots&\dots &\dots &\dots&\dots &\dots\\
q_{M-2}^1&q_{M-2}^2&q_{M-2}^3 &\dots  &I_{n_{M-2}}&0&0\\
q_{M-1}^1&q_{M-1}^2&q_{M-1}^3 &\dots &q_{M-1}^{M-2} &I_{n_{M-1}}&0\\
q_{M}^1&q_{M}^2&q_{M}^3 &\dots &q_M^{M-2}&q_M^{M-1} &I_{n_{M}}
\end{array}
\right).
$$
Its diagonal is formed by the set of $M$ square blocks $n_k\times n_k$. Each diagonal block is
proportional to the unit matrix of the corresponding dimension. Block $q_i^j$ is $n_i\times n_j$ matrix.

Let $[Q]_k$ be its diagonal lower $k\times k$ block
$$[Q]_k:=\left(
\begin{array}{ccccccc}
I_{n_{M-k+1}} &0 &0&\dots &0&0&0\\
q_{M-k+2}^{M-k+1}&I_{n_{M-k+2}} &0 &\dots &0&0&0\\
q_{M-k+3}^{M-k+1}&q_{M-k+3}^{M-k+2}&I_{n_{M-k+3}} &\dots &0&0&0\\
\dots&\dots&\dots &\dots &\dots&\dots &\dots\\
q_{M-2}^{M-k+1}&q_{M-2}^{M-k+2}&q_{M-2}^{M-k+3} &\dots  &I_{n_{-2}}&0&0\\
q_{M-1}^{M-k+1}&q_{M-1}^{M-k+2}&q_{M-1}^{M-k+3} &\dots &q_{M-1}^{M-2} &I_{n_{M-1}}&0\\
q_{M}^{M-k+1}&q_{M}^{M-k+2}&q_{M}^{M-k+3} &\dots &q_M^{M-2}&q_M^{M-1} &I_{n_{M}}
\end{array}
\right),
$$
so $Q=[Q]_M, [Q]_1=I_{n_{M}}$.

Denote the non-trivial parts of the vector-column-blocks by $\vec q_k$:
$$
\vec{q\ }^k:=(q_{k+1}^k, q_{k+2}^k, \dots , q_{M-2}^k, q_{M-1}^k, q_M^k)^T.
$$
They are rectangular matrices of the dimension 
$(n_{k+1}+ n_{k+2}+ \dots + n_{M-2} + n_{M-1} + n_M)\times n_k$.
Consider the vector-raw-blocks $\vec p_k$
$$
\vec p_k:=(p_k^{k+1}, p_k^{k+2}, \dots , p_k^{M-2}, p_k^{M-1}, p_k^M)
$$
 They are $n_k\times (n_{k+1}+ n_{k+2}+ \dots + n_{M-2} + n_{M-1} + n_M)$ matrices. 
 The blocks $q_i^j$ and $p_j^i$
 have the dimensions $n_i\times n_j$ and $n_j\times n_i$ correspondingly.

Consider upper-triangular matrix $\rho$
\begin{equation}\label{ro}
\rho:=\left(
\begin{array}{ccccccc}
\lambda'_1I_{n_1} &\rho_1^2 &\rho_1^3 &\dots &\rho_1^{M-2} &\rho_1^{M-1} &\rho_1^M \\
0&\lambda'_2I_{n_2} &\rho_2^3  &\dots &\rho_2^{M-2} &\rho_2^{M-1} &\rho_2^{M} \\
0&0&\lambda'_3I_{n_3} &\dots &\rho_3^{M-2} &\rho_3^{M-1} &\rho_3^{M} \\
\dots&\dots&\dots &\dots &\dots&\dots &\dots\\
0&0&0&\dots&\lambda'_{n_{M-2}}I_{n_{M-2}}&\rho_{M-2}^{M-1}&\rho_{M-2}^{M}\\
0&0&0&\dots&0&\lambda'_{n_{M-1}}I_{n_{M-1}}&\rho_{M-1}^{M}\\
0&0&0&\dots&0&0&\lambda'_{n_M}I_{n_{M}}
\end{array}
\right).
\end{equation}
Denote non-trivial parts of the vector-raw-blocks by $\vec \rho_k$:
$$
\vec \rho_k:=(\rho_k^{k+1}, \rho_k^{k+2}, \dots , \rho_k^{M-2}, \rho_k^{M-1}, \rho_k^M)
$$

\begin{theorem}\label{thm-formula}
Matrix $A$:
\begin{equation}\label{formula}
A=Q\rho Q^{-1},
\end{equation}
where the block-vector-raws of $\rho$ are 
$$\vec \rho_k:=\vec p_k [Q]_{M-k}
$$
provides the canonical parameterization of the orbit $\mathcal O (J)\ni A$ by the couples of matrix
elements of blocks $p_j^i, q_i^j$: $(p_j^i)_{st}, (q_i^j)_{ts}, \ 1\le s\le n_j, 1\le t\le n_i, \ \ 1\le n_i,n_j \le M  $.
\end{theorem}

\underline{Proof}

To find the Jordan structure of $\rho$ we construct the hierarchy (\ref{eq10}) for it.   

Let us proof that on the open set of the matrix elements of $\rho$ 
the hierarchy (\ref{eq10}) gives the lower-diagonal blocks of $\rho$.

The first $n_1$ columns of $\rho-\lambda_1 I$ vanish. Consider the equality corresponding to (\ref{eq23}) for this stair-flight:
$$
Pg=\hat gJ'+\hat P\tilde J_{\rho}.
$$
Here $\tilde J_{\rho}$ is the normal Jordan form of the lower diagonal block of $\rho$ and $J'$ complements
$\tilde J_{\rho}$ to the normal Jordan form of $\rho$.  

We do not know the $\tilde J_{\rho}$ and $J'$ now, $J'$ may have too many zero columns. 
Consider the columns corresponding to the zero columns of $\tilde J_{\rho}$.
Denote matrices collected from these columns only by $[\![\dots ]\!]$. It is the projection on the subset of columns:
$$[\![ Pg]\!]=[\![\hat gJ']\!].
$$
Matrix $g$ is non-degenerate, consequently on the open set of matrix elements of $P$ matrix $[\![ Pg]\!]$ has a full
rank that is $m-\dim\ker\tilde J_{\rho}$. It implies $\rank J'=m-\dim\ker\tilde J_{\rho}$ or
$$\rank 
\left(\begin{array}{c}J'\\\tilde J_{\rho}\end{array}
\right)=m.
$$
Geometrically it means the following.
\begin{itemize}
\item No one of the Jordan chains of the lower block of $\rho$ in question was finished on the flight 
of the hierarchy. 
\item We started new $n_k-\dim\ker \tilde J_{\rho}=n_k-n_{k+1}$ chains.
\end{itemize} 

It proves that the Jordan structures of $\rho$ and $J$ coincide.  
\vskip 2ex

Let us construct the canonical coordinates for (\ref{formula}) using the method developed in the present paper. 
We proved that the kernel of $\rho-\lambda'_1 I$ 
is formed by the first $n_1$ columns on the open set of matrix elements of $\rho$.

It is easy to verify that
$$Q=\left(
\begin{array}{cc}
I&0\\
\vec{q\ }^1&[Q]_{M-1}
\end{array}
\right)=
\left(
\begin{array}{cc}
I&0\\
\vec{q\ }^1&I
\end{array}
\right)
\left(
\begin{array}{cc}
I&0\\
0&[Q]_{M-1}
\end{array}
\right),
$$
consequently
$$Q^{-1}=
\left(
\begin{array}{cc}
I&0\\
0&([Q]_{M-1})^{-1}
\end{array}
\right)
\left(
\begin{array}{cc}
I&0\\
-\vec{q\ }^1&I
\end{array}
\right).
$$
The substitution of these representations of $Q$ and $Q^{-1}$ to (\ref{formula}) gives the first flight 
of the hierarchy (\ref{eq10}). The diagonal lower block has the same structure as (\ref{formula}).
The iteration of the procedure gives the statement of the theorem.
\qed
\vskip 2ex

The inverse map
$$\mathcal O(J) \to \mathcal E(F_1\oplus E_1)\times \mathcal E(F_2\oplus E_2)\times \dots  \mathcal E(F_M\oplus E_M)
$$
involves the construction of the hierarchy (\ref{eq10}). It is a sequence of the couples of steps. We should find the eigenspace
of the diagonal lower block and change the first part of the basic vectors to the normalized basis of the eigenspace.

\section{Examples}

Let us consider examples. The canonical parameterization of $A\in \mathcal O'$ is given by the product $Q\rho Q^{-1}$.
\begin{example}
Let $N=4$, $\lambda_i=\lambda_j\Leftrightarrow i=j$.
\end{example}
$$Q=\left(
\begin{array}{cccc}
1&0&0&0\\
q_4&1&0&0\\
q_5&q_2&1&0\\
q_6&q_3&q_1&1\\
\end{array}
\right)
$$
$$
Q^{-1}=\left(\left(
\begin{array}{cccc}
1&0&0&0\\
0&1&0&0\\
0&0&1&0\\
0&0&-q_1&1
\end{array}
\right)
\left(
\begin{array}{cccc}
1&0&0&0\\
0&1&0&0\\
0&-q_2&1&0\\
0&-q_3&0&1
\end{array}
\right)\right)
\left(
\begin{array}{cccc}
1&0&0&0\\
-q_4&1&0&0\\
-q_5&0&1&0\\
-q_6&0&0&1
\end{array}
\right)
$$
$$=\left(
\begin{array}{llll}
1&0&0&0\\
-q_4&1&0&0\\
-q_5+q_4q_2&-q_2&1&0\\
-q_6+q_5q_1-q_4(-q_3+q_1q_2)&-q_3+q_1q_2&-q_1&1\end{array}
\right)
$$
The corresponding matrix $\rho$ is
$$\left(
\begin{array}{cccc}
\lambda_4&p_4+p_5q_2+p_{6}q_3&p_5+p_{6}q_1&p_{6}\\
0&\lambda_3&p_2+p_3q_1&p_3\\
0&0&\lambda_2&p_1\\
0&0&0&\lambda_1\\
\end{array}
\right)
$$

\begin{example}
Let $N=5$, $\lambda_i=\lambda_j\Leftrightarrow i=j$
\end{example}
$$Q=\left(
\begin{array}{ccccc}
1&0&0&0&0\\
q_7&1&0&0&0\\
q_8&q_4&1&0&0\\
q_9&q_5&q_2&1&0\\
q_{10}&q_6&q_3&q_1&1\\
\end{array}
\right)
$$
The corresponding matrix $\rho$ is
$$\left(
\begin{array}{ccccc}
\lambda_5&p_7+p_8q_4+p_9q_5+p_{10}q_6&p_8+p_9q_2+p_{10}q_3&p_9+p_{10}q_1&p_{10}\\
0&\lambda_4&p_4+p_5q_2+p_{6}q_3&p_5+p_{6}q_1&p_{6}\\
0&0&\lambda_3&p_2+p_3q_1&p_3\\
0&0&0&\lambda_2&p_1\\
0&0&0&0&\lambda_1\\
\end{array}
\right)
$$
\begin{example}
The Jordan box $4\times 4$ with zero eigenvalue
\end{example}
$$Q=\left(
\begin{array}{cccc}
1&0&0&0\\
q_4&1&0&0\\
q_5&q_2&1&0\\
q_6&q_3&q_1&1\\
\end{array}
\right), \rho=
\left(
\begin{array}{cccc}
0&p_4+p_5q_2+p_{6}q_3&p_5+p_{6}q_1&p_{6}\\
0&0&p_2+p_3q_1&p_3\\
0&0&0&p_1\\
0&0&0&0\\
\end{array}
\right)
$$

\begin{example}
Let $N=6$, 
$J=
\left(
\begin{array}{cccccc}
0&1&0&0&0&0\\
0&0&1&0&0&0\\
0&0&0&0&0&0\\
0&0&0&0&1&0\\
0&0&0&0&0&0\\
0&0&0&0&0&1
\end{array}
\right)
$
\end{example}
$$Q=\left(
\begin{array}{cccccc}
1&0&0&0&0&0\\
0&1&0&0&0&0\\
q_{10}&q_6&1&0&0&0\\
q_{11}&q_{7}&0&1&0&0\\
q_{12}&q_8&q_4&q_2&1&0\\
q_{13}&q_9&q_5&q_3&q_1&1
\end{array}
\right)
$$
the corresponding matrix $\rho$ is
$$\left(
\begin{array}{cccccc}
0&0&p_{10}+p_{12}q_4+p_{13}q_5&p_{11}+p_{12}q_2+p_{13}q_3&p_{12}+p_{13}q_1&p_{13}\\
0&0&p_6+p_8q_4+p_{9}q_5&p_7+p_{8}q_2+p_{9}q_3&p_{8}+p_9q_1&p_9\\
0&0&0&0&p_4+p_5q_1&p_5\\
0&0&0&0&p_2+p_3q_1&p_3\\
0&0&0&0&0&p_1\\
0&0&0&0&0&1
\end{array}
\right)
$$

 \subsection{Acknowledgements}  
 The author is grateful to Yuri~I.~Manin for consultations and warm-hearted 
 attention to this work.
 \vskip 3ex
 
 The author thanks Yuli~Rudyak who has drawn his attention to the book \cite{Dusa} 
 and Eugene~Lerman for the discussions.
 
 This research was carried out while the author was a visitor at Max-Planck Institute f\"ur Mathematik Bonn.

\end{document}